\documentclass[12pt]{amsart}
\usepackage{graphicx}
\usepackage{amssymb}
\usepackage{amsmath}
\usepackage{amsthm,amsfonts,bbm}
\usepackage{amscd}
\usepackage{geometry}
\usepackage{epsfig,epstopdf}
\usepackage[backref]{hyperref}
\usepackage{caption}
\usepackage{subcaption}

\newtheorem{thm}{Theorem}[section]
\newtheorem{cor}[thm]{Corollary}
\newtheorem{lem}[thm]{Lemma}
\newtheorem{prob}{Problem}

\newtheorem{conj}{Conjecture}
\theoremstyle{definition}
\newtheorem{defi}[thm]{Definition}
\theoremstyle{remark}

\newtheorem{exm}[thm]{\bf Example}
\numberwithin{equation}{section}
\numberwithin{figure}{section}
\def \b1{\mathbf{1}}
\def \deg{\text{deg}}
\def \diag{\text{diag}}
\def \e{\epsilon}
\def \Q{\mathbb{Q}}
\def \R{\mathbb{R}}
\def \T{\mathsf{T}}

\begin{document}
\title{Degree-similar graphs and cospectral graphs}

\author[Y.-Z. Fan]{Yi-Zheng Fan*}
\address{Center for Pure Mathematics, School of Mathematical Sciences, Anhui University, Hefei 230601, P. R. China}
\email{fanyz@ahu.edu.cn}
\thanks{*Supported by National Natural Science Foundation of China (Grant No. 12331012).}

\author[R.-J. Xing]{Ruo-Jie Xing}
\address{School of Mathematical Sciences, Anhui University, Hefei 230601, P. R. China}
\email{xingrj@stu.ahu.edu.cn}

\author[Y.-L. Zhang]{Yi-Liu Zhang}
\address{School of Mathematical Sciences, Anhui University, Hefei 230601, P. R. China}
\email{zhangyl@stu.ahu.edu.cn}

\author[W. Wang]{Wei Wang$^\sharp$}
\address{School of Mathematics and Statistics, Xi'an Jiaotong University, Xi'an 710049, P. R. China}
\email{wang\_weiw@xjtu.edu.cn}
\thanks{$^\sharp$Corresponding author. Supported by National Natural Science Foundation of China (Grant No. 12371357).}



\begin{abstract}
Let $G$ be a graph with adjacency matrix $A(G)$ and degree matrix $D(G)$, and let $L_\mu(G):=A(G)-\mu D(G)$.
Two graphs $G_1$ and $G_2$ are called \emph{degree-similar} if there exists an invertible matrix $M$ such that $M^{-1} A(G_1) M =A(G_2)$ and $M^{-1} D(G_1) M =D(G_2)$.
In this paper, we address three problems concerning degree-similar graphs proposed by Godsil and Sun.
First, we present a new characterization of degree-similar graphs using degree partition, from which we derive methods and examples for constructing cospectral graphs and degree-similar graphs.
Second, we construct infinite pairs of non-degree-similar trees $G_1$ and $G_2$ such that $tI- L_\mu(G_1)$ and $tI-L_\mu(G_2)$ have the same Smith normal form over $\Q(\mu)[t]$, which provides a negative answer to a problem posed by Godsil and Sun.
Third, we establish several invariants of degree-similar graphs and obtain results on unicyclic graphs that are degree-similar determined.
Lastly we prove that for a strongly regular graph $G$ and any two edges $e$ and $f$ of $G$, $G \backslash e$ and $G \backslash f$ have identical $\mu$-polynomial, i.e., $\det(tI-L_\mu(G \backslash e))=\det(tI-L_\mu(G \backslash f))$, which enables the construction of pairs of non-isomorphic graphs with same $\mu$-polynomial, where $G \backslash e$ denotes the graph obtained from $G$ by deleting the edge $e$.
\end{abstract}

\maketitle

\section{Introduction}
Let $G=(V(G),E(G))$ be a graph with vertex set $V(G)$ and edge set $E(G)$, and let $A(G)$ and $D(G)$ be respectively the adjacency matrix and degree matrix of $G$.
Godsil and Sun~\cite{GodSun24} introduced the notion of degree similar graphs.
Two graphs $G_1$ and $G_2$ are called \emph{degree-similar} if there exists an invertible matrix $M$ such that
\begin{equation}\label{ds}
M^{-1} A(G_1) M = A(G_2), M^{-1} D(G_1) M = D(G_2).
\end{equation}
Clearly, if $G_1$ and $G_2$ are degree-similar, then their adjacency matrices $A$, Laplacian matrices $L:=D-A$, signless Laplacians $Q:=D+A$, and normalized Laplacians $N:=D^{-1/2}AD^{-1/2}$ are all similar, and hence cospectral with respect to the above matrices.
As noted in \cite{GodSun24} or \cite{WangLLX11}, if $G_1$ and $G_2$ are degree-similar and have no isolated vertices, then their Ihara zeta functions are equal. For more on Ihara zeta functions, see \cite{Ter10}.

Degree-similar graphs have a stronger condition than some earlier versions of cospectral graphs.
The \emph{generalized $\alpha$-adjacency matrix} of a graph $G$ is defined to be
$$A_\alpha(G):=A(G)+\alpha J,$$ where $\alpha \in \R$ and $J$ is an all-one matrix.
The \emph{generalized $\alpha$-characteristic polynomial} (or \emph{$\alpha$-polynomial} for short) of $G$ is defined to be
$$\phi(G,t,\alpha)=\det(tI-A_\alpha(G)).$$
Here, we use the prefix `$\alpha$-' to distinguish it from another generalized adjacency matrix or polynomial to be introduced below.
Johnson and Newman proved the following interesting theorem (see \cite{vanDamHae03, vanDamHaeKoo07}).
For further details on the generalized $\alpha$-adjacency matrix or the generalized $\alpha$-characteristic polynomial, refer to \cite{Tut79,JohnNew80,vanDamHae03,HaeSpe04,vanDamHaeKoo07}.

\begin{thm}\label{GenS}
The following statements are equivalent.

{\em (1)} Two graphs $G_1$ and $G_2$ are cospectral with respect to generalized $\alpha$-adjacency matrix for all $\alpha$.

{\em (2)} $A_\alpha(G_1)$ and $A_\alpha(G_2)$ are cospectral for two distinct values of $\alpha$.

{\em (3)} $G_1$ and $Q_2$ are cospectral with respect to the adjacency matrix, and so are their complements.

{\em (4)} There exists an orthogonal matrix $Q$ such that $Q^\top A(G_1) Q=A(G_2)$ and $Q \b1 =\b1$, where $\b1$ denotes the all-one vector.
\end{thm}

Wang and Xu \cite{WangXuEJC06} called the union of the spectrum of $A(G)$ of a graph $G$ and the spectrum of $A(G^c)$ the \emph{generalized spectrum} of $G$, where $G^c$ denotes the complement of the graph $G$.
Wang and his coauthors investigated the problem of graphs determined by generalized spectrum (or equivalently, determined by $\alpha$-polynomial) in a series of papers \cite{WangXuEJC06,WangXuLAA06,Wang13,Wang17} by using walk-matrices and Smith normal forms over the ring of integers.

The \emph{generalized $\mu$-adjacency matrix} of a graph $G$ is defined by
$$L_\mu(G):=A(G)-\mu D(G),$$
 and the \emph{generalized $\mu$-characteristic polynomial} (or \emph{$\mu$-polynomial} for short) of $G$ is defined by Wang et al. \cite{WangLLX11} as follows:
$$\psi(G,t,\mu):=\det(tI-L_\mu(G)).$$
If $G_1$ and $G_2$ have the same $\mu$-polynomial,
then they are cospectral with respect to the adjacency matrix, the Laplacian matrix, the signless Laplacian matrix and the normalized Laplacian matrix.
Wang et al.~\cite{WangLLX11} proved that if $G_1$ and $G_2$ have the same $\mu$-polynomial, then they have the same degree sequence.
The authors also constructed two non-isomorphic degree-similar graphs which are surely cospectral graphs with respect to generalized $\mu$-adjacency matrix for all $\mu$.
There is no similar result for generalized $\mu$-adjacency matrices as Theorem \ref{GenS} for generalized $\alpha$-adjacency matrices.
For example, there exist two cospectral graphs with respect to $A$ and $L$ but not with respect to $Q$ (\cite[Fig. 4]{vanDamHae03}),
also two cospectral graphs with respect to $A$ and $Q$ but not with respect to $D$ (namely they have different degree sequences) (\cite[Table 4, third pair]{GodMc76}).

By Lemma 4.4 of \cite{GodSun24} (see Lemma \ref{DSC}), if $G_1$ and $G_2$ are degree-similar, and one of them is connected, then their complements are also degree similar.
So, in this case, $G_1$ and $G_2$ have the same generalized spectra, and hence have the same $\alpha$-polynomials by Theorem \ref{GenS}, which are called \emph{$A_\alpha$-cospectral}.

In general, if $G_1$ and $G_2$ are degree-similar over $\R$, surely $L_\mu(G_1)$ and  $L_\mu(G_1)$ are similar over $\R(\mu)$, the latter of which is equivalent to that $tI-L_\mu(G_1)$ and $tI-L_\mu(G_2)$ have the same Smith normal forms (abbreviated as SNFs) over $\R(\mu)[t]$.
By Lemma 9.2 of \cite{GodSun24}, $L_\mu(G_1)$ and  $L_\mu(G_1)$ are similar over $\R(\mu)$ if and only if they are similar over $\Q(\mu)$, which implies that
$tI-L_\mu(G_1)$ and $tI-L_\mu(G_2)$ have the same SNF over $\Q(\mu)[t]$ if $G_1$ and $G_2$ are degree-similar.
Clearly, if $tI-L_\mu(G_1)$ and $tI-L_\mu(G_2)$ have the same SNF over $\Q(\mu)[t]$, then $G_1$ and $G_2$ have the same $\mu$-polynomials by considering the last invariant divisors, which are called \emph{$L_\mu$-cospectral}.

\begin{figure}[htbp]
  \centering
    \includegraphics[scale=1.0]{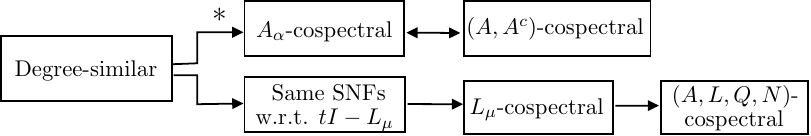}\\
  \caption{\small The implication relations among degree-similarity, SNFs and different versions of cospectral properties} \label{frm}
\end{figure}

By the above discussion, we have the following implication relations listed in Fig. \ref{frm}, where the implication under * means an additional condition of `connectedness', and $(A,A^c)$-cospectral means cospectral with respect to the adjacency matrix $A$ of a graph and the adjacency matrix of the complement of the graph, and $(A,L,Q,N)$-cospectral means cospectral with respect to the djacency matrix $A$, the Laplacian $L$, the signless Laplacian $Q$, and the normalized Laplacian $N$.

Wang et al. \cite{WangLLX11} proposed the following problem:
 \emph{Suppose that two graphs $G_1$ and $G_2$ have the same $\mu$-polynomials, i.e., they are $L_\mu$-cospectral. Does there exist
 an orthogonal matrix $Q$ such that}
\begin{equation}\label{QS}
Q^\top A(G_1) Q=A(G_2), Q^\top D(G_1) Q=D(G_2)?
\end{equation}
Godsil and Sun \cite{GodSun24} give an example of infinite pairs of graphs that share the common $\mu$-polynomials but are not degree-similar, which gives a negative answer to the above problem.
In the same paper \cite{GodSun24}, Godsil and Sun presented three interesting problems on degree-similar graphs as follow.

\begin{prob}\label{P1}\cite{GodSun24}
Find more degree-similar graphs. In particular, are there non-isomorphic
degree-similar unicyclic graphs?
\end{prob}

We give a new characterization of degree-similar graphs by using degree partition, from which we derive some methods for constructing new pairs degree-similar graphs from known ones.
It is known that two trees are degree-similar if and only if they are isomorphic.
Therefore, unicyclic graphs are the first candidates for finding non-isomorphic degree-similar graphs.
By using {\scshape SageMath}, we  could not find non-isomorphic degree-similar unicyclic graphs with at most $20$ vertices.
A graph $G$ is called \emph{degree-similar determined} if any graph that is degree-similar to $G$ must be isomorphic to $G$.
We give some invariants for degree-similar graphs, and prove some classes of unicyclic graphs are degree-similar determined.

\begin{prob}\label{P2}\cite{GodSun24}
Let $G_1$ and $G_2$ be two graphs such that $tI-L_\mu(G_1)$ and $tI-L_\mu(G_2)$ have the same SNF over $\mathbb{Q}(\mu)[t]$.
Are $G_1$ and $G_2$ are degree similar?
\end{prob}

Godsil and Sun \cite{GodSun24} show that if $tI-L_\mu(G_1)$ and $tI-L_\mu(G_2)$ have the same SNF over $\mathbb{Q}(\mu)[t]$, then $A(G_1)$ and $A(G_2)$ are similar over $\mathbb{Q}$, as are $D(G_1)$ and $D(G_2)$.
We give a negative answer to the Problem \ref{P2} by constructing an infinite family of tree pairs.

For a graph $G$ and an edge $e$ of $G$, denote by $G \backslash e$ the graph obtained from $G$ by deleting the edge $e$.
In \cite{GodSunZ23} the authors proved that if $G$ is a strongly regular graph, then
for any two edges $e$ and $f$ of $G$, the graphs $G \backslash e$ and $G \backslash f$ are $(A,L,Q,N)$-cospectral.
In this paper, we prove that $G \backslash e$ and $G \backslash f$ have the same $\mu$- polynomials, or they are $L_\mu$-cospectal, which generalizes Godsil-Sun's result and pushes Problem \ref{P3} a step forward if the answer to Problem \ref{P3} is positive.

\begin{prob} \label{P3}
Let $G$ be a strongly regular graph with two different edges $e$ and $f$. Are
$G \backslash e$ and $G \backslash f$ degree similar?
\end{prob}

The paper is organized as follows.
In Section 2, we present a new characterization of degree-similar graphs by using degree partition, from which we derive some methods and examples for constructing degree-similar graphs or cospectral graphs.
In Section 3, we construct an infinite pairs of non-degree-similar trees $G_1$ and $G_2$ such that $tI-L_\mu(G_1)$ and $tI-L_\mu(G_2)$ share the same SNF, and hence give an negative answer to Problem \ref{P2}.
In Section 4, we give some invariants for degree-similar graphs, and prove some classes of unicyclic graphs are degree-similar determined, pushing the study of Problem \ref{P1}
In Section 5, we prove that for a strongly regular graph $G$ and any two edges $e,f$ of $G$, $G \backslash e$ and $G \backslash f$ have the same $\mu$-polynomials, or they are $L_\mu$-cospectal, pushing the study of Problem \ref{P3}.
Finally we introduce orthogonally degree-similar graphs and some remarks for the notion.

\section{Degree partitions}
In this section we will use degree partition to give a new characterization of degree-similar graphs, from which we present some methods for constructing cospectral graphs and degree similar graphs.
We also give some examples of constructions at the end of this section.

We first introduce some concepts and notations.
Let $G$ be a graph with vertex set $V(G)$, and let $u \in V(G)$.
We use $N_G(u)$ denote the set of neighbors of $u$ in $G$.
The \emph{degree} of $u$, denoted by $\deg_G(u)$, is defined to be the cardinality of the set $N_G(u)$.
Suppose that $G$ has $t$ distinct degrees $d_1,\ldots,d_t$.
The \emph{degree partition} of $G$, denoted by $\pi(G)$, is a partition of the vertex set $V(G)$ of $G$, which consists of subsets
$$V_i=\{v \in V(G): \deg_{G}(v)=d_i\}$$
 for $i \in [t]:=\{1,\ldots, t\}$, namely, $\pi(G)=\{V_1,\ldots,V_t\}$.

Let $M$ be a matrix with rows and columns indexed by the vertices of $G$.
Let $U_1,U_2$ be the subsets of $V(G)$.
Denote by $M[U_1|U_2]$ the submatrix of $M$ with rows indexed by $U_1$ and columns indexed by $U_2$, and by $M(U_1|U_2)$ the submatrix of $M$ with rows indexed by $V(G) \backslash U_1$ and columns indexed by $V(G) \backslash U_2$.
We simply write $M[U_1|U_1]$ as $M[U_1]$ and $M(U_1|U_1)$ as $M(U_1)$.

By Lemma 4.1 of \cite{GodSun24}, the invertible matrix  $M$ in Definition \ref{ds} of degree-similar graphs is block diagonal.
Here we give a more detailed statement by using degree partition.

\begin{lem}\label{DegPart}
Let $G_1,G_2$ be two graphs with same vertex set.
Then $G_1$ and $G_2$ are degree similar if and only if,
by reordering the vertices of $G_1$ and $G_2$, $G_1$ and $G_2$ have the same degree partition, say $\pi=\{V_1,\ldots,V_t\}$, and there exist invertible matrices $M_1,\ldots,M_t$ with rows and columns indexed by $V_1,\ldots,V_t$ respectively, such that
\begin{equation}\label{DSE}
 M_i^{-1} A(G_1)_{ij} M_j = A(G_2)_{ij}, ~ i,j=1,\ldots,t,
\end{equation}
where
\begin{equation}\label{Aij}
A(G_k)_{ij}:=A(G_k)[V_i| V_j], ~ k=1,2;  i,j=1,\ldots, t.
\end{equation}
\end{lem}

\begin{proof}
Suppose that $G_1,G_2$ are degree-similar graphs with the same vertex set $V$.
Then there exists an invertible matrix $M$ such that
$$ M^{-1} D(G_1) M =D(G_2), M^{-1} A(G_1) M =A(G_2).$$
Let $d_1,\ldots,d_t$ be all distinct degrees of $G_1$, and let $U_i=\{v \in V(G_1): \deg_{G_1}(v)=d_i\}$ for $i \in [t]$.
Since  $M^{-1} D(G_1) M=D(G_2)$, the graph $G_2$ has the same degree sequences as $G_1$.
Let $W_i=\{v \in V(G_2): \deg_{G_2}(v)=d_i\}$ for $i \in [t]$.
Note that $|U_i|=|W_i|$ for $i\in [t]$.

Now, by reordering the vertices of $G_1$, for some permutation matrix $P$, we have
$$ P D(G_1) P^\top =\diag\left(d_1 I_{|U_1|}, \cdots, d_t I_{|U_t|} \right)=:\tilde{D}.$$
Similarly, by reordering the vertices of $G_2$, for some permutation matrix $P'$,
$$ P' D(G_2) {P'}^\top =\diag\left(d_1 I_{|W_1|}, \cdots, d_t I_{|W_t|}\right)=\tilde{D}.$$
So, after the above reordering of the vertices, $G_1$ and $G_2$ have the same degree partition, say $\pi=\{V_1,\ldots,V_t\}$.
The matrices $\tilde{A}(G_1)=PA(G_1)P^\top$ and $\tilde{A}(G_2)=P'A(G_2)P'^\top$ are respectively
the adjacency matrices of $G_1$ and $G_2$ after the reordering of vertices.

Let $\tilde{M}=PMP'^\top$.
We have
  $$ \tilde{M}^{-1} \tilde{D} \tilde{M} = \tilde{D},  ~
  \tilde{M}^{-1} \tilde{A}(G_1) \tilde{M} = \tilde{A}(G_2).$$
Partition $\tilde{M}$ conformable with $\pi$, and let $\tilde{M}_{ij}:=\tilde{M}[V_i|V_j]$ for $i,j \in [t]$.
Since $\tilde{D} \tilde{M} = \tilde{D}\tilde{M}$, we have
$$ d_i \tilde{M}_{ij} =   \tilde{M}_{ij} d_j, ~ i,j \in [t],$$
which implies that $\tilde{M}_{ij}=0$ for $i \ne j$.
Hence $\tilde{M}=\diag\{\tilde{M}_{ii}: i \in [t]\}$, a block diagonal compatible with $\pi$.
Let $\tilde{A}(G_k)_{ij}=\tilde{A}(G_k)[V_i| V_j], k=1,2, i,j=1,\ldots, t$.
From the fact $\tilde{M}^{-1} \tilde{A}(G_1) \tilde{M} =\tilde{A}(G_1)$, we have
$$ \tilde{M}_{ii}^{-1} \tilde{A}(G_1)_{ij} \tilde{M}_{jj}= \tilde{A}(G_2)_{ij}, ~ i,j \in [t].$$
The necessity now follows by taking $M_i=\tilde{M}_{ii}$ for $i \in [t]$ and noting $\tilde{A}(G_k)$ is the adjacency matrix of $G_k$ after reordering of vertices for $k =1,2$.

Conversely, if $G_1$ and $G_2$ have the same degree partition $\pi$, then by reordering the vertices, we can write the degree matrices $D(G_1)$ and $D(G_2)$ in the following form:
$$ D(G_1)=D(G_2)=\diag \left(d_1 I_{|V_1|}, \cdots, d_t I_{|V_t|}\right).$$
Let $M=\diag \left(M_1, \ldots, M_t \right)$.
It is easy to verify that
$$ M^{-1} D(G_1) M = D(G_2), ~ M^{-1} A(G_1) M  = A(G_2).$$
So, $G_1$ and $G_2$ are degree-similar.
\end{proof}

\begin{lem}\cite{GodSun24} \label{DSC}
If $G_1,G_2$ are degree-similar and one of them is connected, then their complements are degree-similar.
\end{lem}

In Lemma \ref{DegPart}, if replacing all $M_i$'s by $kM_i$'s for any nonzero $k$, or equivalently replacing $M$ by $kM$, the Eq. \eqref{DSE} still holds, where $M=\diag\{M_i: i \in [t]\}$.
If, in addition, one of $G_1$ and $G_2$ is connected, from the proof of Lemma \ref{DSC},
the matrix $M$ satisfies
$$ M^{-1} J M = J.$$
So $M$ has constant row sum and constant column sum, implying that
$$ MJ=JM=c J$$
for some nonzero $c$.
By taking $c=1$ we have the following result.

\begin{cor}\label{DSEc}
Let $G_1$ and $G_2$ be two graphs on the same vertex set, where $G_1$ is connected.
Then $G_1$ and $G_2$ are degree-similar if and only if,
by reordering the vertices of $G_1$ and $G_2$, $G_1$ and $G_2$ have the same degree partition, say $\pi=\{V_1,\ldots,V_t\}$, and there exist  invertible matrices $M_1,\ldots,M_t$ with rows and columns indexed by $V_1,\ldots,V_t$ respectively, such that
\begin{equation}\label{DSE2}
 M_i^{-1} A(G_1)_{ij} M_j = A(G_2)_{ij}, M_i^\top \b1=M_i \b1=\b1, i,j=1,\ldots,t,
\end{equation}
where $A(G_k)_{ij}$ is defined as in \eqref{Aij}.
\end{cor}

We give the following result for construction of degree-similar graphs from a known pair of degree-similar graphs.

\begin{thm}\label{DScos}
Let $G_1,G_2$ be degree-similar graphs on the same vertex set, which have the same degree partition $\pi=\{V_1, \ldots, V_t\}$, where $G_1$ is connected.
For $k=1,2$, let $G_k[V_i]$ be the subgraph of $G_k$ induced by $V_i$ for $i \in [t]$, and let $G_k[V_i,V_j]$ be the bipartite subgraph of $G_k$ with vertex sets $V_i \cup V_j$ whose edges are those of $G_k$ connecting $V_i$ and $V_j$ for $i \ne j$ and $i,j \in [t]$.
Let $\tilde{G}_1,\tilde{G}_2$ be obtained from $G_1,G_2$ respectively by applying some of following operations simultaneously:

{\em (1)} replacing some $G_k[V_i]$'s with their complements,

{\em (2)} replacing some $G_k[V_i]$'s with empty graphs,

{\em (3)} replacing some $G_k[V_i,V_j]$'s for $i \ne j$ with their complements in complete bipartite graph with bipartition $\{V_i,V_j\}$,

{\em (4)} replacing some $G_k[V_i,V_j]$'s with empty graphs,

\noindent
Then, with respect to adjacency matrix, $\tilde{G}_1$ is cospectral with $\tilde{G}_2$ with cospectral complements.

Furthermore, if both $\tilde{G}_1$ and $\tilde{G}_2$ have the same degree partition as $\pi$, then $\tilde{G}_1$ is degree similar to $\tilde{G}_2$.
In particular, taking operation (1) if each vertex of $V_i$ has degree $(|V_i|-1)/2$ in the graph $G_k[V_i]$, and (or) taking operation (3) if each vertex of $V_i$ has degree $|V_j|/2$ and each vertex of $V_j$ has degree  $|V_i|/2$ in the graph $G_k[V_i,V_j]$,
then $\tilde{G}_1$ is degree-similar to $\tilde{G}_2$.
\end{thm}

\begin{proof}
By Corollary \ref{DSEc}, there exist invertible matrices $M_i$ with rows and columns indexed by $V_i$ for $i\in [t]$, such that
$$ M_i^{-1} A(G_1)_{ij} M_j = A(G_2)_{ij}, M_i^\top \b1=M_i \b1=\b1, i,j=1,\ldots,t,$$
where $(A_k)_{ij}$ is defined as in \eqref{Aij}.

Let $A(\tilde{G}_k)_{ij}:=A(\tilde{G}_k)[V_i|V_j]$ for $k=1,2$ and $i,j=1,2,\ldots,t$.
Let $M=\diag\{M_i: i \in [t]\}$.
To verify $\tilde{G}_1$ is cospectral with $\tilde{G}_2$, it suffices to prove
$M^{-1} A(\tilde{G}_1) M = A(\tilde{G}_2)$, or equivalently,
\begin{equation}\label{tldM}
M_i^{-1} A(\tilde{G}_1)_{ij} M_j = A(\tilde{G}_2)_{ij}, i,j=1,\ldots,t.
\end{equation}
Observe that
if taking operation (1), $A(\tilde{G}_k)_{ii}=J-I-A(G_k)_{ii}$;
and if taking operation (3), $A(\tilde{G}_k)_{ij}=J-A(G_k)_{ij}$.
Since $M_i^\top \b1=M_i \b1=\b1$, we have
$$ M_i^{-1} A(\tilde{G}_1)_{ii} M_i= M_i^{-1} (J-I-A(G_1)_{ii}) M_i = J-I-A(G_2)_{ii}=A(\tilde{G}_2)_{ii},$$
and
$$ M_i^{-1} A(\tilde{G}_1)_{ij} M_j= M_i^{-1}(J-A(G_k)_{ij})M_j=J-A(G_k)_{ij}=A(\tilde{G}_2)_{ij}.$$
Similarly, if taking operation (2), $A(\tilde{G}_k)_{ii}=O$;
and if taking operation (3), $A(\tilde{G}_k)_{ij}=O$, where $O$ denotes a zero matrix of appropriate size.
Obviously, $M_i^{-1} O M_i = O$, $M_i^{-1} O M_j = O$.
So, Eq. \eqref{tldM} holds, and $\tilde{G}_1$ is cospectral with $\tilde{G}_2$.
Using the fact $M^\top \b1=M \b1=\b1$ and noting $A(G^c)=J-I-A(G)$ for a graph $G$, we have
$$ M^{-1} A (\tilde{G}_1^c) M = A(\tilde{G}_2^c),$$
implying that $\tilde{G}_1$ and  $\tilde{G}_2$ have cospectral complements.

If $\tilde{G}_1$ and  $\tilde{G}_2$ have the same degree partition as $\pi$,
surely $M^{-1} D(\tilde{G}_1) M = D(\tilde{G}_2)$.
Combing with the proved equality \eqref{tldM}, we get $\tilde{G}_1$ is degree-similar to $\tilde{G}_2$.
If each vertex of $V_i$ has degree $(|V_i|-1)/2$ in the graph $G_k[V_i]$,
taking the operation (1) in $G_k$ will preserve the degree of each vertex of $G_k$.
Similarly, if each vertex of $V_i$ has degree $|V_j|/2$ and each vertex of $V_j$ has degree  $|V_i|/2$ in the graph $G_k[V_i,V_j]$, taking operation (3) in $G_k$ also preserves the degree of each vertex of $G_k$.
So, $\tilde{G}_1$ and  $\tilde{G}_2$ have the same degree partition as $\pi$, and hence they are degree-similar.
\end{proof}

By Theorem \ref{DScos}, we will produce $3^t \cdot 3^{t \choose 2}$ pairs of cospectral graphs from a pair of degree-similar graphs $G_1$ and $G_2$, where $t$ is the number of parts in the degree partition of $G_1$ or $G_2$.
Maybe some of these pairs of graphs are isomorphic.
Next we give some examples of cospectral graphs and degree-similar graphs by using Theorem \ref{DScos}.

\begin{exm}
The first pair of non-isomorphic degree-similar graphs $X_{1,1}$ and $X_{1,2}$ in Fig. \ref{DSCOS1} were introduced by Wang et al.~\cite{WangLLX11}.
We use three kinds of colored vertices for degree partition, and denote by $V_r, V_g,V_b$ the set of red vertices of degree $4$, the set of green vertices of degree $3$ and the set of blue vertices of degree $2$, respectively.

\begin{figure}[htbp]
  \centering
    \includegraphics[scale=.8]{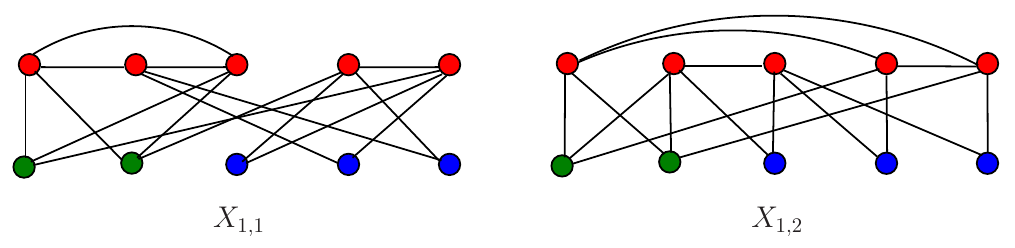}\\
  \caption{\small Degree-similar graphs $X_{1,1}$ and $X_{1,2}$ (\cite{WangXuLAA06})} \label{DSCOS1}
\end{figure}

By taking the complements of $X_{1,k}[V_r]$, we get a pair of cospectral graphs $X_{2,k}$ with cospectral complements for $k=1,2$; see Fig. \ref{DSCOS2}.

\begin{figure}[htbp]
  \centering
    \includegraphics[scale=.9]{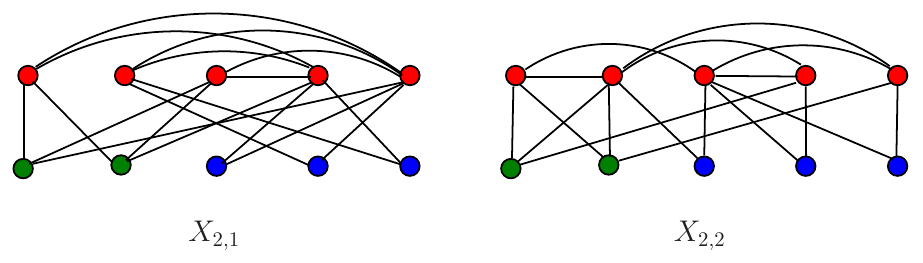}\\
  \caption{\small Cospectral graphs $X_{2,1}$ and $X_{2,2}$ with cospectral complements} \label{DSCOS2}
\end{figure}

By replacing $X_{1,k}[V_r]$ with empty graphs, we get a pair of cospectral graphs $X_{3,k}$ with cospectral complements for $k=1,2$; see Fig. \ref{DSCOS3}.

\begin{figure}[htbp]
  \centering
    \includegraphics[scale=.9]{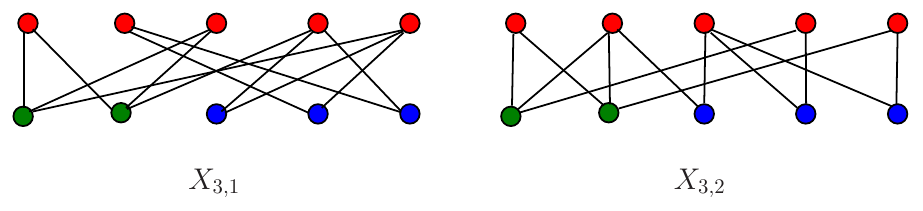}\\
  \caption{\small Cospectral graphs $X_{3,1}$ and $X_{3,2}$ with cospectral complements} \label{DSCOS3}
\end{figure}

By taking complements of $X_{1,k}[V_r,V_g]$ in the complete bipartite graph with two parts $V_r$ and $V_g$, we get a pair of cospectral graphs $X_{4,k}$ with cospectral complements for $k=1,2$; see Fig. \ref{DSCOS4}.

\begin{figure}[htbp]
  \centering
    \includegraphics[scale=.9]{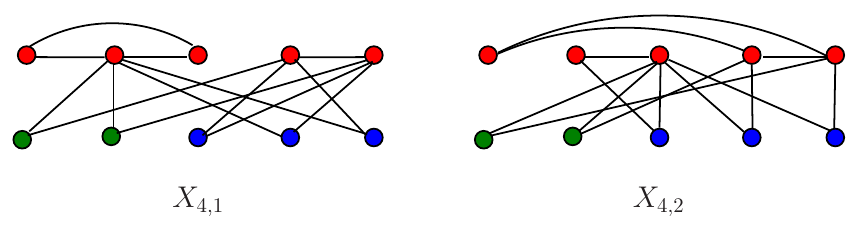}\\
  \caption{\small Cospectral graphs $X_{4,1}$ and $X_{4,2}$ with cospectral complements} \label{DSCOS4}
\end{figure}

If replacing $X_{1,k}[V_r,V_g]$ by empty graphs, we get a pair of cospectral graphs $X_{5,k}$ with cospectral complements for $k=1,2$; see Fig. \ref{DSCOS5}.
By deleting the isolated green vertices, we have two cospectral tricyclic graphs which are isomorphic.

\begin{figure}[htbp]
  \centering
    \includegraphics[scale=.9]{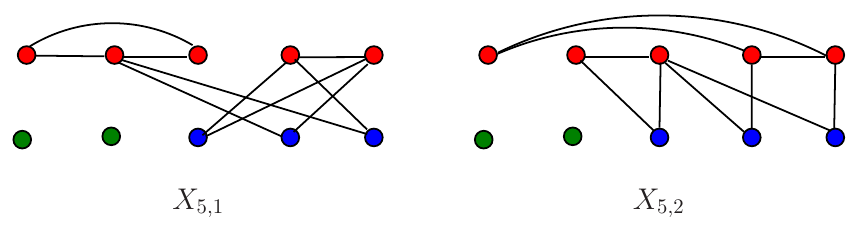}\\
  \caption{\small Cospectral graphs $X_{5,1}$ and $X_{5,2}$ with cospectral complements} \label{DSCOS5}
\end{figure}

If replacing  $X_{4,k}[V_r,V_b]$ by empty graphs, we will get two cospectral graphs $X_{6,k}$ with cospectral complements; see Fig. \ref{DSCOS6}.
By deleting the blue vertices, we get two non-isomorphic cospectral bicyclic graphs.

\begin{figure}[htbp]
  \centering
    \includegraphics[scale=.9]{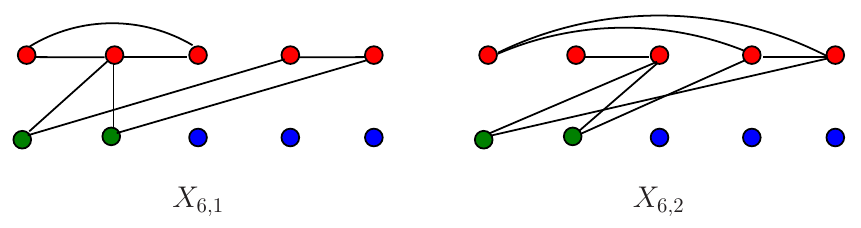}\\
  \caption{\small Cospectral graphs $X_{6,1}$ and $X_{6,2}$ with cospectral complements} \label{DSCOS6}
\end{figure}

\end{exm}

\begin{exm}
The second pair of non-isomorphic degree-similar graphs $Y_{1,1}$ and $Y_{1,2}$ in Fig. \ref{DSGS1} were introduced by Godsil and Sun \cite{GodSun24}.

\vspace{2mm}

\begin{figure}[htbp]
  \centering
    \includegraphics[scale=.9]{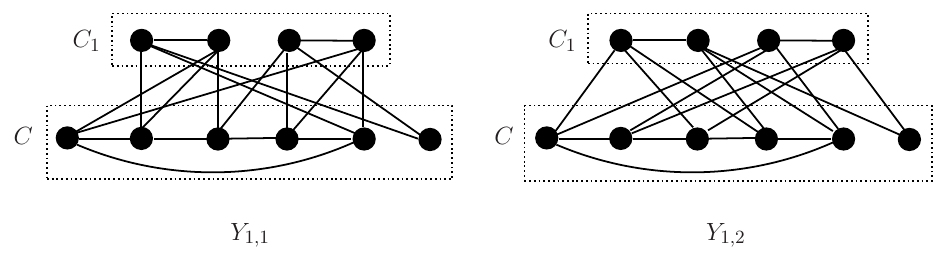}\\
  \caption{\small Degree-similar graphs $Y_{1,1}$ and $Y_{1,2}$ (\cite{GodSun24})} \label{DSGS1}
\end{figure}

By Lemma 6.2 of \cite{GodSun24}, for any graph $Y$, adding all possible edges between $Y$ and $C_1$ (or $Y$ and $C$, or $Y$ and $C\cup C_1$) in Fig. \ref{DSGS1}, the resulting two graphs are also degree-similar.
If letting $Y=P_3$, we get a pair of degree similar graphs $Y_{2,1}$ and $Y_{2,2}$ in Fig. \ref{DSGS2}; see Example 6.3 of \cite{GodSun24}.

\begin{figure}[htbp]
  \centering
    \includegraphics[scale=.8]{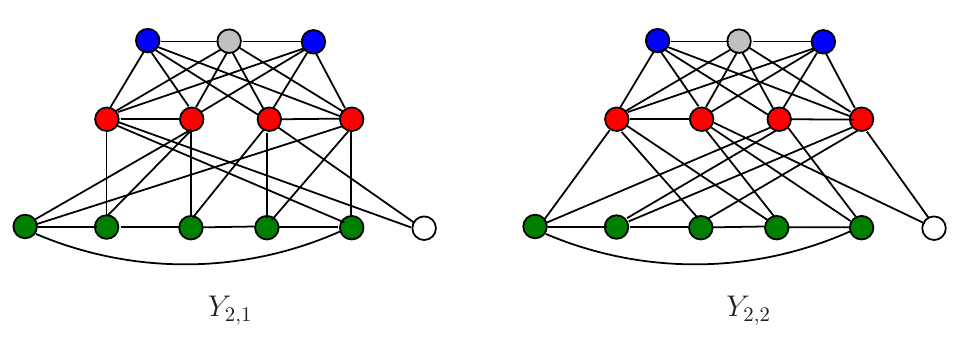}\\
  \caption{\small Degree-similar graphs $Y_{2,1}$ and $Y_{2,2}$ (\cite{GodSun24})} \label{DSGS2}
\end{figure}

By taking complements of the subgraphs induced on green vertices, we get a pair of degree-similar graphs $Y_{3,1}$ and $Y_{3,2}$ by Theorem \ref{DScos}; see Fig. \ref{DSGS3}.
In fact, replacing the path $P_3$ in $Y_{3,k}$ (for $k=1,2$) by any nontrivial connected graph $Y$ and adding all possible edges between $Y$ and $C_1$ (the red vertices), the resulting two graphs are still degree-similar.

\begin{figure}[htbp]
  \centering
    \includegraphics[scale=.8]{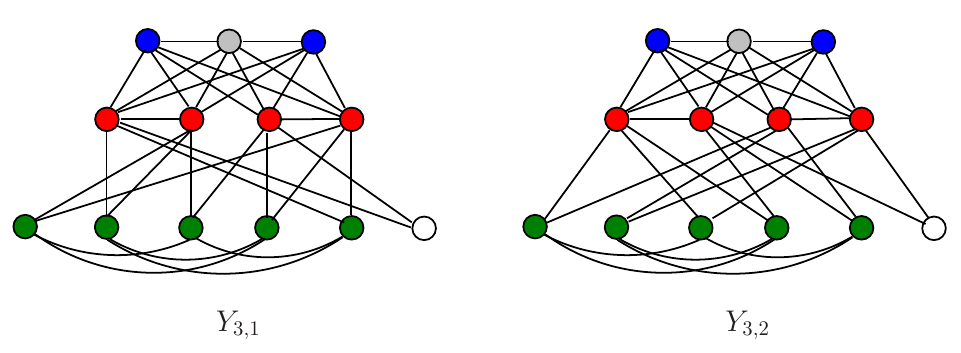}\\
  \caption{\small Degree-similar graphs $Y_{3,1}$ and $Y_{3,2}$} \label{DSGS3}
\end{figure}
\end{exm}

\section{Trees}
In this section, we will construct an infinite family of tree pairs $G_1$ and $G_2$ such that $tI-L_\mu(G_1)$ and $tI-L_\mu(G_2)$ have the same Smith normal form over $\Q(\mu)[t]$ but $G_1$ and $G_2$ are not degree-similar, and hence give a negative answer to Problem \ref{P2} asked by Godsil and Sun \cite{GodSun24}.

\begin{figure}[htbp]
  \centering
    \includegraphics[scale=.8]{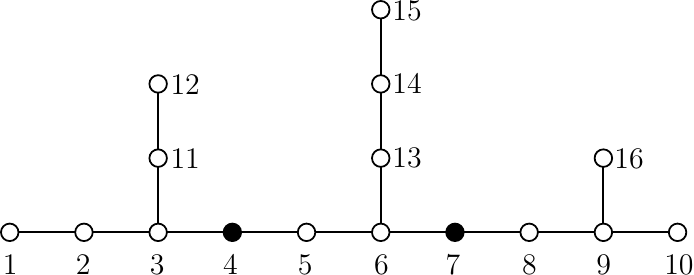}\\
  \caption{\small A tree $\T$ on $16$ vertices} \label{fig_tree}
\end{figure}

Let $G_1$ and $G_2$ be two graphs with roots $u$ and $v$ respectively. The \emph{coalescence} of $G_1$ and $G_2$, denoted by $G_1(u) \odot G_2(v)$, is the graph formed by identifying the root $u$ of $G_1$ and the root $v$ of $G_2$.
The following tree $\T$ in Fig. \ref{fig_tree} was appeared in \cite{McKay77} for constructing non-isomorphism cospectral graphs.
McKay \cite{McKay77} showed that for any nontrivial tree $T$ with root $r$, $T(r) \odot \T(4)$ is not isomorphic to $T(r) \odot \T(7)$, but they are cospectral with respect to adjacency matrix, Laplacian matrix and signless Laplacian matrix, and also normalized Laplacian matrix \cite{Osb13}.

Let $G$ be a general nontrivial graph with root $r$.
Let $G_1:=G(r) \odot \T(4)$ and $G_2:=G(r) \odot \T(7)$.
Godsil and Sun \cite{GodSun24} proved that $G_1$ and $G_2$ have the same $\mu$-polynomial, namely,
$\psi(G_1,t,\mu)=\psi(G_2,t,\mu)$, but $G_1$ is not degree-similar to $G_2$ when $G$ is  any nontrivial tree, which answered a problem proposed by Wang et al.~ \cite{WangLLX11}.
In the following we will prove that $tI- L_\mu(G_1)$ and $tI-L_\mu(G_2)$ have the same Smith normal form when $G$ is a path.

\begin{lem}\label{APath}
let $\T$ be the tree in Fig. \ref{fig_tree}, and let $P_{m+1}$ be a path on $m+1$ vertices with one endpoint $r$, where $m \ge 0$.
Let $G_1:=P_{m+1}(r) \odot \T(4)$ and $G_2:=P_{m+1}(r) \odot \T(7)$.
Then $tI- L_\mu(G_1)$ and $tI-L_\mu(G_2)$ have the same Smith normal form over $\Q(\mu)[t]$.
\end{lem}

\begin{proof}
Let $n:=m+16$ be the number of vertices of $G_1$ or $G_2$.
Along the path $P_{m+1}$ starting from the root $r$, label other vertices of $P_{m+1}$ as $17,18,\ldots,n$, where $n$ is the another endpoint of $P_{m+1}$.
Denote $d_{k,l}(G_i):=\det (tI - L_\mu(G_i)(k|l))$ and $d_{k}(G_i):=\det (tI - L_\mu(G_i)(k))$ for $i=1,2$.

We first investigate the $(n-1)$th determinant divisor of $tI - L_\mu(G_1)$, denoted by $D_{n-1}(G_1)$.
By a direct calculation,
$$ d_{n,17}(G_1)=\det (tI - L_\mu(G_1)[V(\T)]), ~ d_{n,4}(G_1)=\det (tI - L_\mu(G_1)[V(\T)\backslash\{4\}]),$$
and
\begin{equation}\label{GCD}
 \gcd(d_{n,17}(G_1),d_{n,4}(G_1))=\alpha(t,\mu)\beta(t,\mu)\gamma(t,\mu),
\end{equation}
where
\begin{equation}\label{abc}
\begin{split}
& \alpha(t,\mu):=(t + \mu), ~  \beta(t,\mu):=  (t^2 + 3\mu t + 2 \mu^2 - 1), \\
& \gamma(t,\mu):= (t^3 + 6\mu t^2 + (11\mu^2 - 3)t + 6\mu^3 - 5\mu).
\end{split}
\end{equation}
So
$$ D_{n-1}(G_1) \mid \alpha(t,\mu)\beta(t,\mu)\gamma(t,\mu).$$

In the following, by Claims 1-3, we will prove that neither of $\alpha(t,\mu)$, $\beta(t,\mu)$ and $\gamma(t,\mu)$ divides $D_{n-1}(G_1)$, which implies that $D_{n-1}(G_1)=1$.

\textbf{Claim 1:} $ \alpha(t,\mu) \nmid D_{n-1}(G_1)$.
Otherwise, $\alpha(t,\mu) \mid d_{10}(G_1)$.
Expanding  $d_{10}(G_1)$ at the vertex $16$, we have
$$ d_{10}(G_1) =(t+\mu) \det (tI - L_\mu(G_1)(10,16))-\det (tI - L_\mu(G_1)(10,16,9)).$$
Noting that $\alpha(t,\mu)=t+\mu$, we have
$$ (t+\mu) \mid \det (tI - L_\mu(G_1)(10,16,9)).$$
Similarly, expanding the above determinant at the vertices $15, 12, 1, n$ successively, if $n \ge 18$,
$$ (t+\mu) \mid \det (tI - L_\mu(G_1)(10,16,9,15,14,12,11,1,2,n,n-1));$$
and if $n =17$,
$$ (t+\mu) \mid \det (tI - L_\mu(G_1)(10,16,9,15,14,12,11,1,2,17,4)).$$
Let
$$ U=
\begin{cases}
\{10,16,9,15,14,12,11,1,2,n,n-1\}, & \text{~if~} n \ge 18,\\
\{10,16,9,15,14,12,11,1,2,17,4\}, & \text{~if~} n =17.
\end{cases}
$$
Now taking $t=-\mu$, we have
$$ \det (-\mu I - L_\mu(G_1)(U))=\det (\mu(D'-I) -A')=0,$$
where $D',A'$ are the principal submatrices of $D(G_1)$ and $A(G_1)$ indexed by the vertices of $V(G_1) \backslash U$, respectively.
As all the vertices of $V(G_1) \backslash U$ have degree greater than $1$,
each diagonal entry of $D'-I$ is positive.
So, for sufficiently large $\mu$, $\mu(D'-I) -A'$ strictly diagonal dominant, and hence $\mu(D'-I) -A'$ is nonsingular, which yields a contradiction.

\textbf{Claim 2:} $\beta(t,\mu) \nmid D_{n-1}(G_1)$.
Otherwise,  $\beta(t,\mu) \mid d_{1}(G_1)$.
Expanding $\det(tI-L_\mu(G_1))$ at the vertex $1$, we have
\begin{equation}\label{C2}
 \det(tI-L_\mu(G_1))=(t+\mu)d_{1}(G_1) -\det(tI-L_\mu(G_1)(1,2)).
 \end{equation}
As $\beta(t,\mu) \mid \det(tI-L_\mu(G_1))$, we have
$$\beta(t,\mu) \mid \det(tI-L_\mu(G_1)(1,2)).$$
Again, expanding $\det(tI-L_\mu(G_1)(1,2))$ at the vertex $3$, we have
\begin{equation}\label{T4-1}
\begin{split}
\det(tI-L_\mu(G_1)(1,2)) &= (t+3\mu) \det(tI-L_\mu(T)[11,12]) \\
& ~~~ \times \det(tI-L_\mu(G_1)(1,2,3,11,12))\\
& ~~~ -(t+\mu)\det(tI-L_\mu(G_1)(1,2,3,11,12))\\
& ~~~ -\det(tI-L_\mu(T)[11,12]) \det(tI-L_\mu(G_1)(1,2,3,11,12,4)).
\end{split}
\end{equation}
As $\beta(t,\mu)=\det(tI-L_\mu(T)[11,12])$ by a direct calculation,
we have $$\beta(t,\mu) \mid \det(tI-L_\mu(G_1)(1,2,3,11,12)).$$
Expanding $\det(tI-L_\mu(G_1)(1,2,3,11,12))$ at the vertex $13$, we have
\begin{equation}\label{T4-2}
\begin{split}
\det(tI-L_\mu(G_1)(1,2,3,11,12)) &= (t+2\mu) \det(tI-L_\mu(T)[14,15]) \\
& ~~~ \times \det(tI-L_\mu(G_1)(1,2,3,11,12,13,14,15))\\
& ~~~ -(t+\mu)\det(tI-L_\mu(G_1)(1,2,3,11,12,13,14,15))\\
& ~~~-\det(tI-L_\mu(T)[14,15]) \\
& ~~~ \times \det(tI-L_\mu(G_1)(1,2,3,11,12,13,14,15,6)).
\end{split}
\end{equation}
As $\beta(t,\mu)=\det(tI-L_\mu(T)[14,15])$ also by a direct calculation,
we have $$\beta(t,\mu) \mid \det(tI-L_\mu(G_1)(1,2,3,11,12,13,14,15)).$$

If taking $\mu=0$, then $\beta(t,0)=t^2-1$ is a factor of $\det(tI-A(G_1)(W))$, where $W:=\{1,2,3,11,12,13,14,15\}$.
So, $1$ is an eigenvalue of $A(G_1)(W)$.
Note that $A(G_1)(W)=A(G_1(W))$, the adjacency matrix of the subgraph $G_1(W)$ which is obtained from $G_1$ by deleting all vertices of $W$ together with their incident edges.
Let $x$ be an eigenvector of $A(G_1(W))$ corresponding to the eigenvalue $1$.
By eigenvector equation, for each vertex $u \in V(G_1) \backslash W$,
\begin{equation}\label{ev}
 x_u=\sum_{v \in N_{G_1(W)}(u)} x_v.
\end{equation}
So, if letting $x_n=a$, then $x_{n-1}=a$ and $x_{n-2}=0$, and
along the path $P$ from the vertex $n$ to the vertex $9$, the values of part vertices of $G_1(W)$ given by $x$ are listed in Fig. \ref{tree_ev}.

Therefore, $x_9$ has one of the following values:
$a,0,-a$.
If $x_9=a$, then $x_{10}=x_{16}=a$ by eigenvector equation, and hence $x_8=x_9-x_{10}-x_{16}=-a$, which yields  a contradiction as a vertex of $P$ with value $a$ can not be adjacent to a vertex of $P$ with value $-a$.
If $x_9=0$, then $x_{10}=x_{16}=0$, and hence $x_8=0$, also yielding a contradiction.
Similarly, if $x_9=-a$, we also get a contradiction as discussed in the case of $x_9=a$.

\begin{figure}[htbp]
  \centering
    \includegraphics[scale=.8]{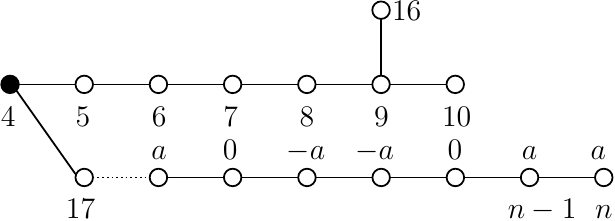}\\
  \caption{\small The graph $G_1(W)$ and part entries of eigenvector $x$} \label{tree_ev}
\end{figure}

\textbf{Claim 3:} $\gamma(t,\mu) \nmid D_{n-1}(G_1)$.
Otherwise, $ \gamma(t,\mu) \mid d_{10,16}(G_1)$.
By a direct calculation,
$$ d_{10,16}(G_1)=\det (tI- L_\mu(G_1))(9,10,16).$$
So, $$ \gamma(t,\mu) \mid \det (tI- L_\mu(G_1))(9,10,16).$$
Expanding $\det (tI - L_\mu(G_1))$ at the vertex $9$,
\begin{equation}\label{CLM3}
\begin{split}
\det (tI - L_\mu(G_1))  & = (t+3\mu) (t+\mu)^2 \det (tI - L_\mu(G_1))(9,10,16)\\
& ~~~ - 2(t+\mu) \det(tI - L_\mu(G_1))(9,10,16) \\
& ~~~ - (t+\mu)^2 \det(tI - L_\mu(G_1))(8,9,10,16),
\end{split}
\end{equation}
which implies
$$ \gamma(t,\mu) \mid  \det(tI - L_\mu(G_1))(8,9,10,16).$$
Again, expanding $\det (tI- L_\mu(G_1))(9,10,16)$ at the vertex $8$, we have
$$ \gamma(t,\mu) \mid \det(tI - L_\mu(G_1))(7,8,9,10,16);$$
and expanding $\det (tI- L_\mu(G_1))(8,9,10,16)$ at the vertex $7$,
$$ \gamma(t,\mu) \mid \det(tI - L_\mu(G_1))(6,7,8,9,10,16).$$
Note that
\begin{align*}
\det(tI - L_\mu(G_1))(6,7,8,9,10,16) & =
\det(tI-L_\mu(T)[13,14,15]) \\
& ~~~ \times  \det(tI-L_\mu(G_1)[\{1,2,3,4,5,11,12\} \cup V(P_{m})]),
\end{align*}
and $\gamma(t,\mu)$ is coprime to $\det(tI-L_\mu(T)[13,14,15])$,
where $P_m$ is a that subpath of $P_{m+1}$ obtained by removing the root $r$.
So
$$ \gamma(t,\mu) \mid \det(tI-L_\mu(G_1)[\{1,2,3,4,5,11,12\} \cup V(P_{m})]).$$
Expanding the above determinant at the vertex $4$, we have
\begin{align*}
& \det(tI-L_\mu(G_1)[\{1,2,3,4,5,11,12\} \cup V(P_{m})])  \\
& = (t+3\mu) (t+2\mu) \det(tI-L_\mu(T)[1,2,3,11,12]) \det (tI-L_\mu(G_1)[V(P_{m-1})])\\
& ~~~ - \det(tI-L_\mu(T)[1,2,3,11,12]) \det (tI-L_\mu(G_1)[V(P_{m})])\\
& ~~~ -(t+2\mu) \det(tI-L_\mu(T)[1,2])^2 \det (tI-L_\mu(G_1)[V(P_{m})]) \\
&~~~ - (t+2\mu) \det(tI-L_\mu(T)[1,2,3,11,12]) \det (tI-L_\mu(G_1)[V(P_{m-1})]),
\end{align*}
where $P_{m-1}$ is the subpath of $P_{m}$ by removing the endpoint $17$.
By a direct calculation,
$\gamma(t,\mu)$ divides $\det(tI-L_\mu(T)[1,2,3,11,12])$ and is coprime to
$(t+2\mu) \det(tI-L_\mu(T)[1,2])^2$.
We have
\begin{equation}\label{C3}
\gamma(t,\mu) \mid \det (tI-L_\mu(G_1)[V(P_{m})]).
\end{equation}
As $\gamma(t,\mu)$ has degree $3$ in $t$, we can assume $m \ge 3$; otherwise we would have a contradiction.

If taking $t=\mu$, then $\gamma(\mu,\mu)=8\mu(3\mu^2-1)$ will divide
$$\delta_m(\mu):=\det (\mu I-L_\mu(G_1)[V(P_{m})])=\det(\mu(D'+I)-A'),$$
 where $D'$ is a degree diagonal matrix on the vertices $P_m$ with entries $D'_{uu}=2$ for all $v \ne n$ and $D'_{nn}=1$, and $A'$ is the adjacency matrix of $P_{m}$.
So
$$ \delta_m(1/\sqrt{3})=\det (1/\sqrt{3}\cdot (D'+I)-A')=0.$$
 This implies that the matrix $1/\sqrt{3}\cdot (D'+I)-A'$ has an eigenvector $x$ corresponding to the zero eigenvalue.
  By eigenvector equation, for all the vertices $u$ of $P_m$ other than $n$,
 $$ \sqrt{3} x_u =\sum_{ v \in N_{P_m}(u)} x_v,$$
 and for the last vertex $n$,
 $$ 2/\sqrt{3} \cdot x_n = x_{n-1}.$$
So, if letting $x_{17}=a$, then $x_{18}=\sqrt{3}a$ and $x_{19}=2a$.
Along the path $P_m$ starting from the vertex $17$, the values of part vertices of $P_{m}$ given by $x$ are listed in Fig. \ref{tree_ev2}.

Therefore, the value $x_{n-1}$ belongs to the set $S:=\{\pm a, \pm \sqrt{3}a, \pm 2a, 0\}$.
It suffices to consider the cases of $x_{n-1}$ having values among $a,\sqrt{3}a,2a,0$.
If $x_{n-1}=a$, then $x_n=\sqrt{3}a/2$, and hence $x_{n-2}=\sqrt{3}a/2$, yielding a contradiction as $x_{n-2} \in S$.
Similarly, if $x_{n-1}=\sqrt{3}a$, then $x_n=x_{n-2}=3a/2$;
and if $x_{n-1}=0$, then $x_n=x_{n-2}=0$ and then $x=0$; which also yields contradiction.
For the last case, if $x_{n-1}=2a$, then $x_n=x_{n-2}=\sqrt{3}a$, and
$$x_{n-3}=a, x_{n-4}=0,x_{n-5}=-a.$$
So, in this case, we have
\begin{equation}\label{evenM}
m \equiv 4 \mod 12.
\end{equation}

\vspace{3mm}

 \begin{figure}[htbp]
  \centering
    \includegraphics[scale=1.0]{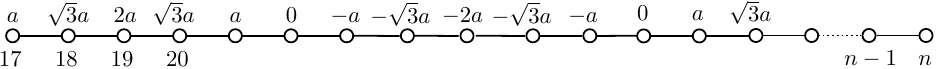}\\
\caption{\small The path $P_m$ and part entries of the eigenvector $x$} \label{tree_ev2}
\end{figure}

If taking $t=-\mu$, then $\gamma(-\mu,\mu)=-2 \mu$ will divide the following determinant
 $$\det (-\mu I-L_\mu(G_1)[V(P_{m})])=\det(\mu(D'-I)-A')=:\eta_m(\mu),$$
  where $D'$ and $A'$ are defined as in the above.
By recursive formula,
$$ \eta_m(\mu) = \mu \eta_{m-1}(\mu) -\eta_{m-2}(\mu).$$
So, if $\mu \mid \eta_m(\mu)$, so does $\eta_{m-2}(\mu)$.
As $m$ is even by Eq. \eqref{evenM}, then $\mu \mid \eta_2(\mu)$.
However, $\eta_2(\mu)=-1$, which yields a contradiction.

Next, along a similar line, by Claims 4-6, we will prove the $(n-1)$th determinant divisor of $tI - L_\mu(G_2)$, denoted by $D_{n-1}(G_2)$, is also $1$.
By a direct calculation,
$$ d_{n,17}(G_2)=\det(tI-L_\mu(G_2)[V(T)])=\det(tI-L_\mu(G_1)[V(T)])=d_{n,17}(G_1), $$
$$d_{n,7}(G_2)=\det(tI-L_\mu(G_2)[V(T)\backslash\{7\}])
=\det(tI-L_\mu(G_1)[V(T)\backslash\{4\}])=d_{n,4}(G_1),$$
and hence by \eqref{GCD},
$$\gcd(d_{n,17}(G_2),d_{n,7}(G_2))=\alpha(t,\mu)\beta(t,\mu)\gamma(t,\mu),$$
where $\alpha(t,\mu),\beta(t,\mu),\gamma(t,\mu)$ are defined as in \eqref{abc}.

{\bf Claim 4:} $\alpha(t,\mu) \nmid D_{n-1}(G_2)$. Otherwise, by expanding $d_{10}(G_2)$ at the vertex $16$, we have
$$\alpha(t,\mu) \mid \det(tI-L_\mu(G_2)(10,16,9)),$$
 and successively expanding determinants at the vertex $15,12,1,n$, if $n \ge 18$,
$$
\alpha(t,\mu) \mid \det(tI-L_\mu(G_2)(10,16,9,15,14,12,11,n,n-1))$$
and if $n=17$,
$$
\alpha(t,\mu) \mid \det(tI-L_\mu(G_2)(10,16,9,15,14,12,11,17,4)).$$
We will get a contradiction by taking $t=-\mu$ and a similar discussion as in the last part of Claim 1.

{\bf Claim 5:} $\beta(t,\mu) \nmid D_{n-1}(G_2)$.
Otherwise, expanding $\det(tI-L_\mu(G_2)$ at the vertex $1$, we have
$$ \beta(t,\mu) \mid \det(tI-L_\mu(G_2)(1,2),$$
expanding $ \det(tI-L_\mu(G_2)(1,2)$ at the vertex $3$, we have
$$ \beta(t,\mu) \mid \det(tI-L_\mu(G_2)(1,2,3,11,12),$$
and expanding $\det(tI-L_\mu(G_2)(1,2,3,11,12)$ at the vertex $13$, we have
$$ \beta(t,\mu) \mid \det(tI-L_\mu(G_2)(1,2,3,11,12,13,14,15).$$

Now taking $\mu=0$, then $\beta(t,0)=t^2-1$ is factor of the determinant
$\det(tI-A(G_2(W))$, where $W=\{1,2,3,11,12,13,14,15\}$, which implies that
the adjacency matrix $A(G_2(W))$ has an eigenvalue $1$.
Let $x$ be an eigenvector of $A(G_2(W))$ corresponding to the eigenvalue $1$.
If $x_{10}:=a$, then $x_9=x_{16}=a$, $x_8=-a$, and $x_7=-2a$.
If $x_4:=b$, then $x_5=b$, $x_6=0$ and $x_7=-b$.
So $2a=b$, and hence $x_{17}=-a$.
The values of $x$ of part vertices of $G_2(W)$ are listed in Fig. \ref{tree_ev4}.
However, $x_{n-1}=x_n$ by eigenvector equation, implying $a=0$ and hence $x=0$; a contradiction.

 \begin{figure}[htbp]
  \centering
    \includegraphics[scale=.8]{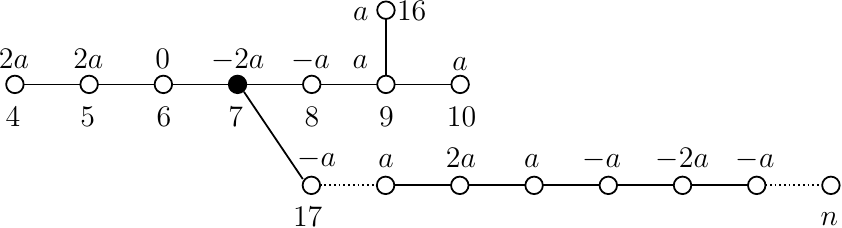}\\
  \caption{\small The graph $G_2(W)$ and part entries of the eigenvector $x$} \label{tree_ev4}
\end{figure}

{\bf Claim 6:} $\gamma(t,\mu) \nmid D_{n-1}(G_2)$.
Otherwise, $\gamma(t,\mu) \mid d_{10,16}(G_2)$.
Note that
 $$d_{10,16}(G_2)=\det (tI-L_\mu(G_2)(9,10,16)),$$
 implying that $\gamma(t,\mu) \mid \det (tI-L_\mu(G_2)(9,10,16))$.
Now, expanding $\det (tI-L_\mu(G_2))$ at the vertex $9$ in a similar way as \eqref{CLM3},
 we have
 $$ \gamma(t,\mu) \mid \det (tI-L_\mu(G_2)(8,9,10,16)).$$
 Expanding $\det (tI-L_\mu(G_2)(8,9,10,16))$ at the vertex $4$, we have
 \begin{align*}
& \det (tI-L_\mu(G_2)(8,9,10,16)) \\
& = (t+2\mu) \det(tI-L_\mu(T)[1,2,3,11,12]) \det (tI - L_\mu(G_2)[\{5,6,7,13,14,15\}\cup V(P_m)])\\
& - \det(tI-L_\mu(T)[1,2])^2 \det (tI - L_\mu(G_2)[\{5,6,7,13,14,15\}\cup V(P_m)])\\
& - \det(tI-L_\mu(T)[1,2,3,11,12]) \det (tI - L_\mu(G_2)[\{6,7,13,14,15\}\cup V(P_m)]).
\end{align*}
As $\gamma(t,\mu)$ divides $\det(tI-L_\mu(T)[1,2,3,11,12])$ and is coprime to $\det(tI-L_\mu(T)[1,2])^2$,
$$ \gamma(t,\mu) \mid \det (tI - L_\mu(G_2)[\{5,6,7,13,14,15\}\cup V(P_m)]).$$
Expanding the above determinant at the vertex $7$, we have
\begin{align*}
& \det (tI - L_\mu(G_2)[\{5,6,7,13,14,15\}\cup V(P_m)]) \\
& ~~~ = (t+3\mu) \det(tI-L_\mu(T)[5,6,13,14,15]) \det (tI - L_\mu(G_2)[V(P_m)])\\
& ~~~ - \det(tI-L_\mu(T)[5,13,14,15]) \det (tI - L_\mu(G_2)[V(P_m)])\\
& ~~~ - \det(tI-L_\mu(T)[5,6,13,14,15]) \det (tI - L_\mu(G_2)[V(P_{m-1})]).
\end{align*}
By a direct calculation,
$\gamma(t,\mu)$ divides $\det(tI-L_\mu(T)[5,6,13,14,15])$ and is coprime to $\det(tI-L_\mu(T)[5,13,14,15])$.
We have
$$ \gamma(t,\mu) \mid \det (tI - L_\mu(G_2)[V(P_m)]),$$
which is consistent with \eqref{C3} in Claim 3.
We will get a contradiction by the same discussion to \eqref{C3}.

By Claims 1-3 and Claims 4-6, we have $D_{n-1}(G_1)=D_{n-1}(G_2)=1$.
By \cite[Lemma 3.1]{GodSun24},
$$ \det (tI-L_\mu(G_1)) = \det (tI-L_\mu(G_2))=:\psi(t,\mu).$$
So $tI-L_\mu(G_1)$ and  $tI-L_\mu(G_2)$  have the same Smith normal form over $\mathbb{Q}(\mu)$ as follows:
$$ 1,\ldots,1, \psi(t,\mu),$$
with $1$ appears $n-1$ times.
So the lemma follows.
\end{proof}

By a similar discussion, we can show  Lemma \ref{APath} also holds if $P_{m+1}$ is replaced by a star with its center as root.
However, due to the length of paper, we omit the result and its proof here.
We believe Lemma \ref{APath} holds when $P_{m+1}$ is replaced by any nontrivial tree.

\begin{conj}
let $\T$ be the tree in Fig. \ref{fig_tree}, and let $T$ be any nontrivial tree with root $r$.
Let $G_1=T(r) \odot \T(4)$ and $G_2=T(r) \odot \T(7)$.
Then $tI- L_\mu(G_1)$ and $tI-L_\mu(G_2)$ have the same Smith normal forms over $\mathbb{Q}(\mu)[t]$.
\end{conj}

We give a negative answer to Problem \ref{P2} asked by Godsil and Sun \cite{GodSun24} by the fact that two trees are degree similar if and only if they are isomorphic \cite{McKay77}.

\begin{cor}
let $\T$ be the tree in Fig. \ref{fig_tree}, and let $P_{m+1}$ be a path on at least $2$ vertices with an endpoint $r$ as root.
Let $G_1=P_{m+1}(r) \odot \T(4)$ and $G_2=P_{m+1}(r) \odot \T(7)$.
Then $tI- L_\mu(G_1)$ and $tI-L_\mu(G_2)$ have the same Smith normal forms over $\mathbb{Q}(\mu)[t]$, but $G_1$ and $G_2$ are not degree similar.
\end{cor}

\section{Unicyclic graphs}
Recall a graph is called \emph{unicyclic} if it is connected and contains only one cycle.
In this section we first give some invariants for degree-similar graphs, and then prove some results on degree-similar determined unicyclic graphs.
If two graphs are degree-similar, then they have same spectra with respect to adjacency matrix, Laplacian matrix, signless Laplacian matrix and normalized Laplaican matrix (if there exist no isolated vertices), respectively.
So we have many invariants for degree-similar graphs, some of which are listed below.

\begin{lem}\label{inv}
Let $G_1$ and $G_2$ be a pair of degree-similar graphs.
Then the following statements hold.

{\em (1)} $G_1$ and $G_2$ have the same numbers of vertices, isolated vertices, edges, connected components, bipartite connected components, respectively.

{\em (2)} If $G_1$ and $G_2$ are connected, then they have the same number of spanning trees.

{\em (3)} If $G_1$ and $G_2$ are connected, then they have the same number of walks of any given length.

\end{lem}

\begin{proof}
By definition, $G_1$ and $G_2$ have the same number of vertices.
Surely, they have the same degree sequence, implying they have the same number of isolated vertices, and also same number of edges as the sum of degrees of a graph is twice the number of edges.
Also by definition, $G_1$ and $G_2$ have the same spectra with respect to Laplacian matrix and signless Laplacian matrix, respectively.
It is well known that the multiplicity of zero as a Laplacian eigenvalue (respectively, as a signless Laplacian eigenvalue) of a graph equals the number of its connected components (respectively, the number of bipartite connected components); see Propositions 1.3.7 and 1.3.9 in \cite{BrouHae11}.
Therefore, $G_1$ and $G_2$ have the same numbers of connected components and bipartite connected components, respectively.

By Matrix-Tree Theorem (or see  Propositions 1.3.4 in \cite{BrouHae11}),
the number of spanning trees of a graph equals the product of nontrivial Laplacian eigenvalues divided by the number of the vertices of the graph.
So, $G_1$ and $G_2$ have the same number of spanning trees if they are connected.

By Corollary \ref{DSEc}, there exists an invertible matrix $M$ such that
$$ M^{-1} A(G_1) M = A(G_2), ~ M^\top \b1 = M \b1 =\b1.$$
Thus, for any positive integer $k$,
$$ \b1^\top A(G_2)^k \b1= \b1^\top M^{-1} A(G_1)^k M \b1=\b1^\top A(G_1)^k \b1,$$
which implies that $G_1$ and $G_2$ have the same number of walks of length $k$.
\end{proof}

Recall that the \emph{girth} of a graph is the minimum length of the cycles in the graph.

\begin{lem}\label{girth}
Let $U_1$ and $U_2$ be two degree-similar unicyclic graphs.
Then they have the same girth.
\end{lem}

\begin{proof}
The result follows by Lemma \ref{inv} (2), since the girth of a unicyclic graph is exactly the number of its spanning trees.
\end{proof}

\begin{thm}
Let $U$ be unicyclic graph on $n$ vertices with girth $g \in \{n, n-1, n-2\}$.
Then $U$ is degree-similar determined, namely, any graph $G$ that is degree-similar to $U$ must be isomorphic to $U$.
\end{thm}

\begin{proof}
Let $G$ be a graph that is degree-similar to $U$.
By Lemma \ref{inv} (1), $G$ is connected with $n$ vertices and $n$ edges, which implies that $G$ is also unicyclic.
By Lemma \ref{girth}, $G$ is a unicyclic graph with the same girth $g$ as $U$.
Thus, if $g=n$ or $g=n-1$, $G$ is isomorphic to $U$ obviously.

Now we consider the case of $g$ equal to $n-2$.
In this case, $U$ has exactly $2$ vertices outside its cycle $C$ of length $n-2$.
Thus, $U$ is one of the following graphs:
$C(r) \odot P_3(u)$, $C(r) \odot P_3(w)$, and $C(r_1,r_2,d)$,
where $P_3$ is a path on $3$ vertices with an endpoint $u$ and a non-endpoint $w$, $C(r_1,r_2,d)$ is obtained from $C$ by attaching one pendent edge at the vertex $r_1$ of $C$ and another pendent edge at $r_2$ of $C$, and the distance between $r_1$ and $r_2$ is $d \ge 1$.
Since $G$ shares the same degree sequence with $U$, if $U=C(r) \odot P_3(u)$ with only one vertex of degree $3$, surely $G \cong U$.
Similarly, if $U=C(r) \odot P_3(w)$ with one vertex of maximum degree $4$, we also have $G \cong U$.

If $U=C(r_1,r_2,d)$, then $G=C(r'_1,r'_2,d')$ for some vertices $r'_1,r'_2$ of $C$ with distance $d'$ by considering the degree sequence.
We assert $d=d'$ and then $G \cong U$.
Otherwise, without loss of generality, assume that $d < d'$.
Let $\omega_{d+2}(U)$ and $\omega_{d+2}(G)$ be the numbers of walks of length $d+2$ in the graph $U$ and $G$, respectively, and let $\omega_{d+2}^{(i)}(U)$ and $\omega_{d+2}^{(i)}(G)$ be the numbers of walks of length $d+2$ in the graph $U$ and $G$ that contain $i$ pendent vertices, respectively, where $i=0,1,2$.
It is easily verified that
$$ \omega_{d+2}^{(0)}(U)=\omega_{d+2}^{(0)}(G), ~ \omega_{d+2}^{(1)}(U)=\omega_{d+2}^{(1)}(G).$$
Observe that the distance between two pendent vertices of $U$ is exactly $d+2$, while the distance between two pendent vertices of $G$ is $d'+2$.
Since $ d < d' \le (n-2)/2$, we have
$$ \omega_{d+2}^{(2)}(U)=1> \omega_{d+2}^{(2)}(G)=0.$$
Therefore,
$$ \omega_{d+2}(U) = \sum_{i=0}^2 \omega_{d+2}^{(i)}(U) > \sum_{i=0}^2 \omega_{d+2}^{(i)}(G) = \omega_{d+2}(G),$$
which yields a contradiction to Lemma \ref{inv} (3).
\end{proof}

Finally in this section we give another class of degree-similar determined unicyclic graphs by using Lemma \ref{DegPart}.

\begin{thm}
Let $T$ be a tree with root $v$, where $T$ contains no vertices of degree $2$ and $v$ is the unique vertex of $T$ with maximum degree.
Let $C_g$ be a cycle of length $g$ with root $r$.
Then the unicyclic graph $U=C_g(r) \odot T(v)$ is degree-similar determined.
\end{thm}

\begin{proof}
Let $G$ is a graph that is degree-similar to $U$.
By Lemma \ref{inv}(1) and Lemma \ref{girth}, $G$ is a unicyclic graph with girth $g$.
By Corollary \ref{DSEc}, we can assume that $G$ and $U$ have the same vertex set, and same degree partition, say $\pi=\{V_1,V_2,\ldots,V_t\}$.
By the assumption on $U$, we can assume $V_1=V(C_g) \backslash \{r\}$, the set of vertices of $U$ with degree $2$; and $V_2=\{r\}$ (or $\{v\}$), the set of the unique vertex of $U$ with maximum degree $2+\deg_T(v)$.
Also, we can write $G=C_g(r) \odot T'(w)$.
By Lemma  \ref{DegPart}, there exist invertible matrices $M_1,\ldots,M_t$ such that
\begin{equation}\label{GenG}
 M_i^{-1} A(U)_{ij} M_j = A(G)_{ij}, i,j \in [t],
\end{equation}
where $A(U)_{ij}=A(U)[V_i|V_j]$ and $A(G)_{ij}=A(G)[V_i|V_j]$ for $i,j \in [t]$.

Observe that $T$ and $T'$ share the same degree partition $\pi'=\{V_2,\ldots,V_t\}$, which is obtained from $\pi$ only by removing $V_1$.
By \eqref{GenG} for $i,j=2,\ldots,t$ and using Lemma \ref{DegPart},
$T$ is degree-similar to $T'$.
By Lemma \ref{inv}(1), $T'$ is also a tree, and hence $T \cong T'$ (\cite{GodMc76}).
As $v$ is unique vertex of $T$ with maximum degree, $w$ is unique vertex of $T'$ with the same maximum degree, and then we have $G \cong U$.
\end{proof}

\section{Strongly regular graphs}
Recall that a graph $G$  is called \emph{strongly regular} with parameters $(n,d;a,c)$ if it has $n$ vertices and is regular of degree $d$, any two adjacency vertices share exactly $a$ common neighbors, and any non-adjacency vertices share exactly $c$ common vertices.
Godsil, Sun and Zhang \cite{GodSunZ23} proved that if $G$ is a strongly regular graph, then
for any two edges $e$ and $f$ of $G$, the graphs $G \backslash e$ and $G \backslash f$ are $(A,L,Q,N)$-cospectral.
In \cite{GodSun24} the authors proposed Problem \ref{P3}, namely, \emph{are
 $G \backslash e$ and $G \backslash f$  degree similar}?
In this section, we will prove that $G \backslash e$ and $G \backslash f$ are $L_\mu$-cospectral, which generalized Godsil-Sun-Zhang's result (\cite[Theorem 1]{GodSunZ23}) and push the  Problem \ref{P3} a step forward.

A graph $G$ is called \emph{walk regular} if for any positive integer $k$, the number of closed walks of length $k$ is the same at all vertices.
If further, the number of walks from vertex $u$ to $v$ of length $k$ is the same for all adjacent vertex pairs $u, v$, then we say $G$ is \emph{$1$-walk regular}.
Surely, a $1$-walk regular graph is regular and  also strongly regular.

Let $G$ be a $1$-walk regular
By Lemma 2.2 of \cite{GodSunZ23}, for any function $f$ defined on the eigenvalues of $A:=A(G)$, there exist $\alpha_f$ and $\beta_f$ such that
\begin{equation}\label{1walk}
 f(A) \circ I = \alpha_f I,  f(A) \circ A = \beta_f A,
\end{equation}
where $\circ$ denotes Schur product.
Let $G$ be a graph and let $u,v$ be vertices of $G$.
Denote by $\delta_{u,v}$ the Kronecker notation, i.e., $\delta_{u,v}=1$ if $u=v$ and  $\delta_{u,v}=0$ otherwise,
and denote by $e_u$ the column vector with rows indexed by the vertices of $G$ whose entries are given by $e_u(v)=\delta_{u,v}$.
Denote $\e_{u,v}=-\mu \delta_{u,v}+ (1- \delta_{u,v})$.

We first give a general result by using a similar technique in \cite{GodSunZ23}.
We need following matrix results for preparation.

\begin{thm}[Sherman-Morrison] \label{SM} Suppose $B$ is an $n \times n$ invertible real matrix and $u,v $ be $n$-dimensional real vectors.
Then $B+uv^\top$ is invertible if and only if $1+v^\top B^{-1} u \ne 0$. In this case,
$$ \left(B+uv^\top \right)^{-1}=B^{-1}-\frac{B^{-1}uv^\top B^{-1}}{1+v^\top B^{-1} u}.$$
\end{thm}

\begin{lem}\cite{GodRoy01}\label{comm}
Assume that $C$ and $D^\top$ are both matrices of size $m \times n$.
Then
$$ \det(I_m-CD) = \det (I_n - DC).$$
\end{lem}

\begin{lem}\label{GCLem}
Let $G$ be a $1$-walk regular graph with adjacency matrix $A$ and degree matrix $D$.
Let $u_1, v_1, \ldots, u_r, v_r$ be vertices in the same clique of $G$.
Then the value of
\begin{equation}\label{DelEdge}
e_{u_r}^\top \left(tI - A + \mu D  \pm (\e_{u_1,v_1} e_{u_1} e_{v_1}^\top + \cdots + \e_{u_{r-1},v_{r-1}} e_{u_{r-1}} e_{v_{r-1}}^\top )\right)^{-1} e_{v_r}
\end{equation}
is independent on the choice of the clique and on the ordering of vertices of the
chosen clique.
\end{lem}

\begin{proof}
Suppose that $G$ is $d$-regular.
Then $ tI - A + \mu D = (t+ \mu d)I -A$.
Let $$f(x)=(t+\mu d - x)^{-1},$$
 which is defined on all eigenvalues of $A$.
As $G$ is $1$-walk regular, by \eqref{1walk}, there exists $\alpha(t,\mu)$ and $\beta(t,\mu)$ such that
$$ f(A) \circ I = \alpha(t,\mu)I,  ~ f(A) \circ A = \beta(t,\mu) A.$$

We prove the result by induction. When $r=1$, as $u_1,v_1$ are in the same clique of $G$,
$$ e_{u_1}^\top f(A) e_{v_1} = \delta_{u_1,v_1} \alpha(t,\mu) + (1-\delta_{u_1,v_1}) \beta(t,\mu),$$
which only depends on whether $u_1$ and $v_1$ are the same or not.

Let $M_0=(t+ \mu d)I -A$ and for $s=1,\ldots,r-1$,
\begin{equation}\label{M}
M_s=t I -A + \mu D \pm (\e_{u_1,v_1} e_{u_1} e_{v_1}^\top + \cdots + \e_{u_s,v_s} e_{u_{s}} e_{v_{s}}^\top).
\end{equation}
Then \eqref{DelEdge} can be written as $e_{u_r}^\top M_{r-1}^{-1} e_{u_r}$.
Assume the result holds for $r=k$, where $k \ge 1$.
Now, by Theorem \ref{SM},
\begin{align*}
e_{u_{k+1}}^\top M_{k}^{-1} e_{v_{k+1}}
& = e_{u_{k+1}}^\top (M_{k-1} \pm \e_{u_k,v_k} e_{u_k} e_{v_k}^\top)^{-1}e_{u_{k+1}} \\
& =  e_{u_{k+1}}^\top \left(M_{k-1}^{-1} \mp \frac{\e_{u_k,v_k} M_{k-1}^{-1}e_{u_k} e_{v_k}^\top M_{k-1}^{-1}}{1+\e_{u_k,v_k}e_{v_k}^\top M_{k-1}^{-1}e_{u_k} }\right)  e_{v_{k+1}} ~~~ \text{(by Theorem \ref{SM})}\\
& = e_{u_{k+1}}^\top M_{k-1}^{-1} e_{v_{k+1}} \mp \frac{\e_{u_k,v_k} (e_{u_{k+1}}^\top M_{k-1}^{-1} e_{u_k}) (e_{v_k}^\top M_{k-1}^{-1} e_{v_{k+1}})}{1+\e_{u_k,v_k}e_{v_k}^\top M_{k-1}^{-1}e_{u_k} },
\end{align*}
whose value does not depend on which clique the vertices are in, and remains
unchanged if we reorder the vertices insides the clique, since each term satisfies
this condition by the induction hypothesis.
\end{proof}

\begin{thm}\label{GCThm}
Let $G$ be a $1$-walk regular graph with clique number $\omega$.
Then for any graph $H$ on at most $\omega$ vertices,
removing edges of $H$ from cliques of $G$ results in graphs with same $\mu$-polynomials.
\end{thm}

\begin{proof}
Let $\bar{H}$ be the graph obtained from $H$ by adding $|V(G)|-|V(H)|$ isolated vertices.
Order the vertices of $\bar{H}$ so that the vertices of $H$ correspond to a clique in $G$.
Assume $H$ has $m$ edges labelled as $e_i=\{u_i,v_i\}$ for $i \in [m]$.
The $\mu$-polynomial of $\hat{G}:=G -E(H)$ is
\begin{equation}
\begin{split}
\psi(\hat{G},t,\mu) & := \det (tI-A+\mu D + (e_{u_1} e_{v_1}^\top + e_{v_1} e_{u_1}^\top \cdots + e_{u_m} e_{v_m}^\top + e_{v_m} e_{u_m}^\top) \\
  & ~~~ ~~~ - \mu (e_{u_1} e_{u_1}^\top + e_{v_1} e_{v_1}^\top + \cdots + e_{u_m} e_{u_m}^\top + e_{v_m} e_{v_m}^\top) ).
\end{split}
\end{equation}
We will prove that $\psi(\hat{G},t,\mu)$ is independent of which clique of $G$ the vertex set of $H$ correspond to or how the vertices of $H$ are ordered.

When $m=1$,
\begin{align*}
\psi(\hat{G},t,\mu) &= \det \left(tI-A+\mu D + (e_{u_1} e_{v_1}^\top + e_{v_1} e_{u_1}^\top) -  \mu (e_{u_1} e_{u_1}^\top + e_{v_1} e_{v_1}^\top)\right)\\
 & =  \det \left(tI-A+\mu D + (e_{u_1}, e_{v_1})
       (e_{v_1} - \mu e_{u_1}, e_{u_1} - \mu e_{v_1} )^\top\right)\\
 &=   \det (tI-A+\mu D) \det \left(I +   (tI-A+\mu D)^{-1}  (e_{u_1}, e_{v_1})
       (e_{v_1} - \mu e_{u_1}, e_{u_1} - \mu e_{v_1} )^\top \right)\\
&=   \det (tI-A+\mu D) \det \left(I_2 +  (e_{v_1} - \mu e_{u_1}, e_{u_1} - \mu e_{v_1} )^\top (tI-A+\mu D)^{-1}  (e_{u_1}, e_{v_1}) \right) \\
 & = \det M_0 \det \left(\begin{array}{cc}
 1+ e_{v_1}^\top M_0^{-1}  e_{u_1} - \mu e_{u_1}^\top M_0^{-1} e_{u_1} &
 e_{v_1}^\top M_0^{-1} e_{v_1} -\mu e_{u_1}^\top M_0^{-1} e_{v_1}\\
 e_{u_1}^\top M_0^{-1} e_{u_1} -\mu e_{v_1}^\top M_0^{-1}  e_{u_1} &
 1+ e_{u_1}^\top M_0^{-1} e_{v_1} - \mu e_{v_1}^\top M_0^{-1}  e_{v_1}
 \end{array} \right),
 \end{align*}
 where the fourth equality follows from Lemma \ref{comm} and $M_0:=tI-A+\mu D$ in the last equality.
In this case, the $\mu$-polynomial
$\psi(\hat{G},t,\mu)$ is independent of the choice of the edge $\{u_1,v_1\}$, since each entry in the $ 2 \times 2$ matrix does not by Lemma \ref{GCLem}.

Define $M_s$ as in \eqref{M}, we have
$$ M_{4m}=tI-A+\mu D + (e_{u_1} e_{v_1}^\top + e_{v_1} e_{u_1}^\top \cdots + e_{u_m} e_{v_m}^\top + e_{v_m} e_{u_m}^\top) - \mu (e_{u_1} e_{u_1}^\top + e_{v_1} e_{v_1}^\top + \cdots + e_{u_m} e_{u_m}^\top + e_{v_m} e_{v_m}^\top).$$
Then
\begin{align*}
\psi(\hat{G},t,\mu)
& = \det \left(M_{4(m-1)}+ (e_{u_m} e_{v_m}^\top + e_{v_m} e_{u_m}^\top) -  \mu (e_{u_m} e_{u_m}^\top + e_{v_m} e_{v_m}^\top)\right)\\
& = \det \left(M_{4(m-1)}+ (e_{u_m}, e_{v_m}) (e_{v_m} - \mu e_{u_m}, e_{u_m} - \mu e_{v_m})^\top \right)\\
 & = \det M_{4(m-1)} \det \left(I + M_{4(m-1)}^{-1} (e_{u_m}, e_{v_m}) (e_{v_m} - \mu e_{u_m}, e_{u_m} - \mu e_{v_m})^\top \right)\\
 &= \det M_{4(m-1)} \det \left(I_2 + (e_{v_m} - \mu e_{u_m}, e_{u_m} - \mu e_{v_m})^\top M_{4(m-1)}^{-1} (e_{u_m}, e_{v_m}) \right)\\
 &= \det M_{4(m-1)}  \\
 & ~~~ \times \det \left(\!\!\! \begin{array}{cc}
 1+ e_{v_m}^\top M_{4(m-1)}^{-1}  e_{u_m} - \mu e_{u_m}^\top M_{4(m-1)}^{-1} e_{u_m} &
 e_{v_m}^\top M_{4(m-1)}^{-1} e_{v_m} -\mu e_{u_m}^\top M_{4(m-1)}^{-1} e_{v_m}\\
 e_{u_m}^\top M_{4(m-1)}^{-1} e_{u_m} -\mu e_{v_m}^\top M_{4(m-1)}^{-1}  e_{u_m} &
 1+ e_{u_m}^\top M_{4(m-1)}^{-1} e_{v_m} - \mu e_{v_m}^\top M_{4(m-1)}^{-1}  e_{v_m}
 \end{array}\!\!\! \right).
 \end{align*}
Since by induction the first factor, and by Lemma \ref{GCLem} each entry in the $ 2 \times 2$ matrix do not dependent on the choice of the clique in $G$ nor on the ordering
of vertices of $H$, the result follows.
\end{proof}

\begin{cor}
Let $G$ be a strongly regular graph with clique number $\omega$ and let $H$ be any graph on at most $\omega$ vertices.
Removing edges of $H$ from cliques of $G$ results in graphs with same $\mu$-polynomials, whose complements also have the same $\mu$-polynomials.
\end{cor}

\begin{proof}
By Theorem \ref{GCThm}, removing edges of $H$ from cliques of $G$ results in graphs with same $\mu$-polynomials.
For the function $f(x)$ defined in Lemma \ref{GCLem},
$$ f(A) \circ I = \alpha(t,\mu) I, f(A) \circ A = \beta(t,\mu) A.$$
Furthermore, if $G$ has parameters $(n,d;a,c)$, then
$A^2 = dI + a A + c(J-I-A)$.
So, there exists $\gamma(t,\mu)$ such that
$$f(A) \circ (J-I-A) = \gamma(t,\mu)  (J-I-A).$$
Let $G^c$ be the complement of $G$, which is also strongly regular or $1$-walk regular.
By a similar arguement in the proof of Theorem \ref{GCThm}, adding edges of $H$ inside a coclique of $G^c$ results in graphs with same $\mu$-polynomials.
Now, deleting edges of $H$ in a clique of $G$ corresponds to adding edges of $H$ in the corresponding coclique of $\bar{G}$.
So, removing edges of $H$ from cliques of $G$ results in graphs whose complements also have the same $\mu$-polynomials.
\end{proof}

\begin{cor}
Let $G$ be a strongly regular graph.
Then for any two edges $e$ and $f$ of $G$, the graphs $G \backslash e$ and $G \backslash f$ have the same $\mu$-polynomial, or they are $L_\mu$-cospectral.
\end{cor}

There are exactly $15$ non-isomorphic strongly regular graphs, denoted by $X_i$ for $i=0,1,\ldots,14$ in \cite{GodSunZ23}, with parameters $(25,12; 5, 6)$.
Their adjacency matrices can be found at Spence's website:
\href{http://www.maths.gla.ac.uk/~es/srgraphs.php}{http://www.maths.gla.ac.uk/\textasciitilde es/srgraphs.php}.
In \cite{GodSunZ23} the authors give a table that lists the number of pairwise non-isomorphic subgraphs of $X_i$ obtained by deleting edges of $6$ small graphs respectively in cliques of $X_i$ for $i=0,1,\ldots,14$.
For example, removing an edge from $X_1$ gives a family of $150$ graphs,
they are pairwise non-isomorphic but $L_\mu$-cospectral.

\section{Orthogonally degree-similar graphs}
Motivated by the problem proposed by Wang et al.~\cite{WangLLX11}, we introduce orthogonally degree-similar graphs, which may be viewed as a stronger version of degree-similar graphs.

\begin{defi}\label{QDS}
Two graphs $G_1$ and $G_2$ are called \emph{orthogonally degree-similar} if there exists an orthogonal matrix $Q$ such that Eq.~\eqref{QS} holds, namely,
$$ Q^\top A(G_1) Q=A(G_2), Q^\top D(G_1) Q=D(G_2).$$
\end{defi}

We have some remarks for the rationality of the definition.

\begin{itemize}

\item[(1)] Two graphs $G_1$ and $G_2$ are cospectral if and only if $A(G_1)$ and $A(G_2)$ are orthogonally similar.

\item[(2)] $G_1$ and $G_2$ are cospectral with cospectral complements if and only if $A(G_1)$ and $A(G_2)$ are similar via an orthogonal matrix $Q$ with $Q \b1=\b1$ (Theorem \ref{GenS}).

\item[(3)] If $G_1$ and $G_2$ are  degree similar and one of them is connected, then $G_1$ and $G_2$ are cospectral with cospectral complements (Lemma \ref{DSC}).
So we have an orthogonal matrix $Q$ as in (2).

\item[(4)] The invertible matrix $M$ in Eq. \eqref{ds} is not unique, since $kM$ still satisfies Eq. \eqref{ds} for any nonzero $k$.

\item[(5)] As the adjacency matrices and degree matrices are symmetric, we have additional requirements for $M$ in Eq. \eqref{ds}, that is,
$$ M^\top A(G_1) (M^{-1})^\top = A(G_2), M^\top D(G_1) (M^{-1})^\top = D(G_2).$$

\item[(6)] The matrices $M$ in most examples of degree similar graphs in \cite{GodSun24} are orthogonal, e.g. Example 5.3, Example 6.3, Example 7.3, Example 8.4.

\end{itemize}

%

\end{document}